\documentclass{amsart}
\usepackage{amsmath,amsthm,amssymb}

\newtheorem{theorem}{Theorem}
\newtheorem{lemma}{Lemma}
\newtheorem{example}{Example}
\newtheorem{problem}{Problem} 
\newtheorem{corollary}{Corollary}
\newtheorem{definition}{Definition}

\newcommand{\bt}{\begin{theorem}}
\newcommand{\et}{\end{theorem}}
\newcommand{\bl}{\begin{lemma}}
\newcommand{\el}{\end{lemma}}
\newcommand{\bex}{\begin{example}}
\newcommand{\eex}{\end{example}}
\newcommand{\bp}{\begin{problem}}
\newcommand{\ep}{\end{problem}}
\newcommand{\bc}{\begin{corollary}}
\newcommand{\ec}{\end{corollary}}
\newcommand{\beq}{\begin{equation}}
\newcommand{\eeq}{\end{equation}}
\newcommand{\benum}{\begin{enumerate}}
\newcommand{\eenum}{\end{enumerate}}
\newcommand{\Z}{\ensuremath{\mathbf Z}}
\newcommand{\N}{\ensuremath{ \mathbf N }}

\title{Problems in Additive Number Theory, I}
\author{Melvyn B. Nathanson}
\address{Department of Mathematics\\
Lehman College (CUNY)\\
Bronx, New York 10468}
\email{melvyn.nathanson@lehman.cuny.edu}
\thanks{Supported in part by grants from the NSA Mathematical Sciences Program and the PSC-CUNY Research Award Program.}
\subjclass[2000]{11B05, 11B13, 11B34}
\keywords{Additive number theory, sumsets, difference sets, representation functions}

\begin{document}

\begin{abstract}
Talk at the Atelier en combinatoire additive (Workshop on Arithmetic Combinatorics) at the Centre de recherches math{\' e}matiques at the Universit{\' e} de Montr{\' e}al  on April 8, 2006.
\end{abstract}

\maketitle

\begin{definition}  
A \emph{problem} is a problem I cannot solve, not necessarily an unsolved problem.
\end{definition}

\section{Sums and differences}
For any set $A$ of integers, we define the \emph{sumset}
\[
A+A =  \{ a + a' : a, a' \in A  \}
\]
and the \emph{difference set}
\[
A-A =  \{ a - a' : a,a' \in A  \}.
\]
In this section we consider finite sets of integers, and the relative sizes of their sumsets and difference sets.
If $A$ is a finite set of integers and $x,y\in \Z$, then the \emph{translation} of $A$ by $x$ is the set  $x+A = \{x+a : a\in A\}$ and the \emph{dilation} 
of $A$ by $y$ is  $y\ast A = \{ya : a\in A\}$.  We have
\[
(x+A) + (x+A) = 2x + 2A 
\]
and
\[
(x+A) - (x+A) = A-A.
\]
Similarly,
\[
y\ast A + y\ast A = y\ast(A+A)
\]
and 
\[
y\ast A - y\ast A = y\ast(A-A).
\]
It follows that 
\[
|(x+y\ast A) + (x+y\ast A)| = |2A|
\]
and
\[
|(x+y\ast A) - (x+y\ast A)| = |A-A|
\]
so the cardinalities of the sum and difference sets of a finite set of integers are invariant under affine transformations of the set.

The set $A$ is \emph{symmetric} with respect to the integer $z$ if $A = z-A$ or, equivalently, if $a\in A$ if and only if $z-a \in A$.  For example, the set $\{4,6,7,9\}$ is symmetric with $z = 13$.  If $A$ is symmetric, then 
\[
A+A = A + (z-A) = z+ (A-A)
\]
and so $|A+A| = |A-A|$.  

If $A = \{a,b,c\}$ with $a<b<c$ and $a+c \neq 2b$, then
\[
|A+A|= 6 < 7 = |A-A|.
\]
If
\[
A = \{ 0,2,3,4,7\}
\]
then 
\begin{align*}
A+A & = [0,14] \setminus \{1,12,13\} \\
A-A & =  [-7,7] \setminus \{-6,6\}
\end{align*}
and
\[
|A+A| = 12 < 13 = |A+A|.
\]

This is the typical situation.  Since
\[
2+7 = 7+2
\]
but
\[
2-7 \neq 7-2
\]
it is natural to expect that in any finite set of integers there are always at least as many differences as sums.
There had been a conjecture, often ascribed incorrectly to 
John Conway\footnote{The confusion may be due to the fact that the first published paper on the conjecture, by John Marica~\cite{mari69}, is entitled "On a conjecture of Conway."  I asked Conway about this at the Logic Conference in Memory of Stanley Tennenbaum at the CUNY Graduate Center on April 7, 2006.  He said that he had actually found a counterexample to the conjecture, and that this is recorded in unpublished notes of Croft~\cite{crof67}.} 
that asserted that  $|A+A| \leq |A-A|$ for every finite set $A$ of integers.  This conjecture is false, and a counterexample is  the set
\beq  \label{pantI:counterexample}
A = \{0,2,3,4,7,11,12,14\},
\eeq
for which
\begin{align*}
A+A &   = [0,28]\setminus \{1,20,27\} \\
A-A &  = [-14,14] \setminus \{\pm 6,\pm 13\}
\end{align*}
and
\[
|A+A| = 26 > 25 = |A-A|.
\]
On the other hand, the conjecture really \emph{should} be true, and suggests the first somewhat philosophical problem.

\bp
Why do there exist finite sets $A$ of integers such that $|A+A| > |A-A|$?
\ep

Given the existence of such aberrant sets, we can ask for the smallest one.  The set $A$ in \eqref{pantI:counterexample} satisfies $|A| = 8$.  

\bp
What is
\[
\min\{ |A| : A \subseteq \Z \text{ and } |A+A| > |A-A| \}?
\]
\ep

Sets $A$ with the property that $ |A+A| > |A-A|$ should have structure.  By structure I do not mean that $A$ contains arithmetic progressions or generalized arithmetic progressions or even subsets of some polynomial type.  There is a significant part of combinatorial and additive number theory, sometimes called \emph{additive combinatorics}, that consists of looking for arithmetic progressions inside sets of integers, or proving that certain sets can be approximated by generalized arithmetic progressions.  The results are beautiful, deep, and difficult, but it is hard to ignore the fact that arithmetic progressions are fundamentally boring, and dense or even relatively dense sets of integers must contain vastly more interesting structures that we have not yet imagined.

The astronomers are trying to understand the large-scale structure of the universe.  If they found an arithmetic progression or a generalized arithmetic progression of galaxies, they would be ecstatic, but it would also be obvious to them that this fascinating and unexpected curiosity is only a small part of the universe, and they would keep looking for other structures.  Since the complexity of sets of integers is comparable to that of the universe, we should also keep looking.

\bp
What is the structure of finite sets satisfying  $|A+A| > |A-A|$?
\ep

If  $A$ is a finite set of integers and $m$ is a sufficiently large positive integer (for example, $m > 2\max(\{ |a| : a\in A\})$, then the set
\[
A_t = \left\{ \sum_{i=0}^{t-1} a_i m^i : a_i \in A \text{ for } i=0,1,\ldots,t-1 \right\}
\]
has the property that 
\[
|A_t+A_t| = |A+A|^t
\]
and
\[
|A_t-A_t| = |A-A|^t.
\]
It follows that if $|A+A| > |A-A|$, then $|A_t+A_t| > |A_t -A_t|$ and, moreover,
\[
\lim_{t\rightarrow\infty} \frac{ |A_t+A_t|}{|A_t-A_t|} = 
\lim_{t\rightarrow\infty} \left(\frac{ |A+A|}{|A-A|}\right)^t = \infty.
\]
The sequence of sets $\{A_t\}_{t=1}^{\infty}$ is the standard parametrized family of sets with more sums than differences.

\bp
Are there other parametrized families of sets satisfying $|A+A| > |A-A|$?
\ep

Even though there exist sets $A$ that have more sums than differences, such sets should be rare, and it must be true with the right way of counting that the vast majority of sets satisfies $|A-A| > |A+A|$.

\bp
Let $f(n)$ denote the number sets $A \subseteq [0,n-1]$ such that $|A-A| < |A+A|$, and let $f(n,k)$ denote the number of such sets $A \subseteq [0,n-1]$ with $|A|=k$.  
Compute 
\[
\lim_{n\rightarrow \infty} \frac{f(n)}{2^n}
\]
and
\[
\lim_{n\rightarrow \infty} \frac{f(n,k)}{{n\choose k}}.
\]
\ep

The functions $f(n)$ and $f(n,k)$ are not necessarily the best functions to count finite sets of nonnegative integers with respect to sums and differences.

\bp
Prove that $|A-A| > |A+A|$ for almost all sets $A$ with respect to other appropriate counting functions.
\ep

\section{Binary linear forms}
The problem of sums and differences can be considered a special case of a more general problem about binary linear forms 
\[
f(x,y) = ux+vy
\]
where $u$ and $v$ are nonzero integers.  For every finite set $A$ of integers, let
\[
f(A) = \{f(a,a') : a,a' \in A\}.
\]
We are interested in the cardinality of the sets $f(A)$.  
For example, the sets associated to the binary linear forms
\[
s(x,y) = x+y
\]
and
\[
d(x,y) = x-y
\]
are the sumset  $s(A) = A+A$ and the difference set $d(A) = A-A$.

We begin by putting our functions into a standard form.
Let $f_0(x,y) = u_0x+v_0y$.  If $d=(u_0,v_0)>1$, let $f_1(x,y) = (u_0/d)x+(v_0/d)y = u_1x+v_1y$.  If $|u_1|<|v_1|$, let $f_2(x,y) = v_1x + u_1y = u_2x+v_2y$.  If $u_2<0$, let $f_3(x,y) = -u_2x-v_2y = u_3x+v_3y$.  
Then
\[
|f_0(A)| =  |f_1(A)| = |f_2(A)| = |f_3(A)| 
\]
for every finite set $A$.   To every binary linear form there is constructed in this way a unique normalized binary linear form $f(x,y) = ux+vy$ such that
\[
u \geq |v| \geq 1 \text{ and } (u,v)=1.
\]
The natural question is:  If $f(x,y)$ and $g(x,y)$ are two distinct normalized binary linear forms, do there exist  finite sets $A$ and $B$ of integers such that $|f(A)| > |g(A)|$ and $|f(B)| < |g(B)|$, and, if so, is there an algorithm to construct $A$ and $B$?  

Brooke Orosz gave constructive solutions to this problem in some important cases.  
For example, she proved the following:  Let  $u > v \geq 1$ and $(u,v)=1,$ and consider the normalized binary linear forms
\[
f(x,y) = ux+vy
\]
and 
\[
g(x,y) = ux-vy.
\]
For $u \geq 3$, the sets
\[
A = \{ 0, u^2-v^2,u^2,u^2+uv\}
\]
and
\[
B = \{0,u^2-uv,u^2-v^2,u^2\}
\]
satisfy the inequality
\[
|f(A)| = 14 > 13 = |g(A)|
\] 
and
\[
|f(B)| = 13 < 14 = |g(B)|.
\] 
For $u =2$, we have $f(x,y) = 2x+y$ and $g(x,y) = 2x-y$.  The sets
\[
A = \{ 0, 3,4,6\}
\]
and
\[
B = \{0,4,6,7\}
\]
satisfy the inequality
\[
|f(A)| = 13> 12 = |g(A)|
\] 
and
\[
|f(B)| = 13 < 14 = |g(B)|.
\] 

The problem of pairs of binary linear forms has been completely solved by Nathanson, O'Bryant, Orosz, Ruzsa, and Silva~\cite{noors06}.

\bt
Let $f(x,y)$ and $g(x,y)$ be distinct normalized binary linear forms.  There exist finite sets $A$, $B$, $C$ with $|C| \geq 2$ such that
\begin{align*}
|f(A)| &  > |g(A)|  \\
|f(B)| &  < |g(B)|  \\
|f(C)| &  = |g(C)|.
\end{align*}
\et

\bp
Let $f(x,y)$ and $g(x,y)$ be distinct normalized binary linear forms. 
Determine if $|f(A)| > |g(A)|$ for most or for almost all finite sets of integers.
\ep

These results should be extended to linear forms in three or more variables.

\bp  
Let $f(x_1,\ldots,x_n) = u_1x_1 + \cdots + u_nx_n$ and $g(x_1,\ldots,x_n) = v_1x_1 + \cdots + v_nx_n$ be linear forms with integer coefficients.  Does there exist a finite set $A$ of integers such that $|f(A)| > |g(A)|$?
\ep

\section{Polynomials over finite sets of integers and congruence classes}

An integer-valued function is a function $f(x_1,x_2,\ldots,x_n)$ such that if $x_1,x_2,\ldots,x_n \in\Z$, then $f(x_1,x_2,\ldots,x_n) \in\Z$.
The binomial polynomial
\[
{x \choose k} = \frac{x(x-1)(x-2)\cdots (x-k+1)}{k!}
\]
is integer-valued, and every integer-valued polynomial is a linear combination with integer coefficients of the polynomials ${x\choose k}$.
For any set $A \subseteq \Z$, we define
\[
f(A) = \{ f(a_1,a_2,\ldots,a_n) : a_i \in A \text{ for } i = 1,2,\ldots,n\} \subseteq \Z.
\]

\bp  \label{pantI:polyineq}
Let $f(x_1,\ldots,x_n)$ and $g(x_1,\ldots,x_n)$ be integer-valued polynomials.  Determine if there exist finite sets $A, B, C$ of positive integers with $|C| \geq 2$ such that 
\begin{align*}
|f(A)| & > |g(A)| \\
|f(B)| & < |g(B)| \\
|f(C)| & = |g(C)| .
\end{align*}
\ep

There is a strong form of Problem~\ref{pantI:polyineq}.
\bp  \label{LF2:prob4}
Let $f(x_1,\ldots,x_n)$ and $g(x_1,\ldots,x_n)$ be integer-valued polynomials.  Does there exist a sequence $\{A_i\}_{i=1}^{\infty}$ of finite sets of integers such that
\[
\lim_{i\rightarrow\infty} \frac{|f(A_i)|}{|g(A_i)|} = \infty ?
\]
\ep

There is also the analogous modular problem.
For every polynomial  $f(x_1,x_2,\ldots,x_n)$ with integer coefficients and for every set $A \subseteq \Z/m\Z$, we define
\[
f(A) = \{ f(a_1,a_2,\ldots,a_n) : a_i \in A \text{ for } i = 1,2,\ldots,n\} \subseteq \Z/m\Z.
\]

\bp  
Let $f(x_1,\ldots,x_n)$ and $g(x_1,\ldots,x_n)$ be polynomials with integer coefficients and let $m \geq 2$.  Do there exist sets $A, B, C \subseteq \Z/m\Z$  with $|C|>1$ such that 
\begin{align*}
|f(A)| & > |g(A)| \\
|f(B)| & < |g(B)| \\
|f(C)| & = |g(C)| .
\end{align*}
\ep

\bp  \label{LF2:prob-modular}
Let $f(x_1,\ldots,x_n)$ and $g(x_1,\ldots,x_n)$ be polynomials with integer coefficients.  Let $M(f,g)$ denote the set of all integers $m \geq 2$ such that there exists a finite set $A$ of congruence classes modulo $m$ such that $|f(A)|  > |g(A)|.$  Compute $M(f,g)$.
\ep

Note that if there exists a finite set $A$ of integers with $|f(A)|  > |g(A)|,$
then $M(f,g)$ contains all sufficiently large integers.

\section{Representation functions of asymptotic bases}
A central topic in additive number theory is the study of bases for the integers and for arbitrary abelian groups and semigroups, written additively.  The set $A$ is called an \emph{additive basis of order $h$} for the set $X$ if every element of $X$ can be written as the sum of exactly $h$ not necessarily distinct elements of $A$.  The set $A$ is called an \emph{asymptotic basis of order $h$} for the set $X$ if all but at most finitely many elements of $X$ can be written as the sum of $h$ not necessarily distinct elements of $A$.  The classical bases in additive number theory for the set $\mathbf{N}_0$ of nonnegative integers are the squares (Lagrange's theorem), the cubes (Wieferich's theorem), the $k$-th powers (Waring's problem and Hilbert's theorem), the polygonal numbers (Cauchy's theorem), and the primes (Shnirel'man's theorem for sufficently large integers).   These classical results in additive number theory are in Nathanson~\cite{nath96,nath00}.

The {\it representation function} for a set $A$ is the function $r_{A,h}(n)$
that counts the number of representations of $n$ as the sum of $h$ elements of $A$.  More precisely, $r_{A,h}(n)$ is the number of $h$-tuples $(a_1, a_2,\ldots, a_h) \in A^h$ such that
\[
a_1 + a_2 + \cdots + a_h = n
\]
and
\[
a_1 \leq a_2 \leq \cdots \leq a_h.
\]
The set $A$ is an asymptotic basis of order $h$ if $r_{A,h}(n) \geq 1$ for all but finitely many elements of $X$.  A fundamental unsolved problem in additive number theory is the \emph{classification problem for representation functions}.

\bp 
Let $f:  \mathbf{N}_0 \rightarrow \mathbf{N}_0 \cup \{\infty\}$ be a function.  Find necessary and sufficient conditions on $f$ in order that there exists a set $A$ in $\mathbf{N}_0$ such that $r_{A,h}(n) = f(n)$ for all $n \in \mathbf{N}_0$.  
\ep

A special case is the classification problem for representation functions for asymptotic bases for the nonnegative integers.

\bp 
Let $\mathcal{F}_0(\N_0)$ denote the set of all functions $f:  \mathbf{N}_0 \rightarrow \mathbf{N}_0 \cup \{\infty\}$ such that $f(n)=0$ for only finitely many nonnegative integers $n$.  For what functions $f \in \mathcal{F}_0(\N_0)$ does there exist a set $A \subseteq \mathbf{N}_0$ such that $r_{A,h}(n) = f(n)$ for all $n \in \mathbf{N}_0$?  
\ep

Nathanson~\cite{nath03} introduced these problems, and recently began to study the representation functions of asymptotic bases for the set $\mathbf{Z}$ of  integers.  He proved~\cite{nath05} that if $f:  \mathbf{Z} \rightarrow \mathbf{N}_0 \cup \{\infty\}$ is {\it any} function such that $f(n)=0$ for only finitely many integers $n$, then there exists a set $A$ in $\mathbf{Z}$ such that
$r_{A,h}(n) = f(n)$ for all integers $n.$  Moreover, arbitrarily sparse infinite sets $A$ can be constructed with the given representation function $f$.  The important new problem is to determine the maximum density of a set $A$ of integers with given representation function $f$. 
For any set $A$ of integers, we define the {\it counting function} 
\[
A(x) = \sum_{a\in A \atop |a| \leq x} 1.
\]

\bp
Let $f:  \mathbf{Z} \rightarrow \mathbf{N}_0 \cup \{\infty\}$ be any function such that $f(n)=0$ for only finitely many integers $n$.
Let $\mathcal{R}(f)$ denote the set of all sets $A \subseteq \Z$ such that
\[
r_{A,2}(n) = f(n)  \text{ for all } n \in \Z.
\]
Compute
\[
\sup\left\{\alpha : A \in  \mathcal{R}(f) 
\text{ and } A(x) \gg x^{\alpha} \text{ for all } x \geq x_0\right\}.
\]
\ep

Nathanson~\cite{nath04} proved that for any function $f$ there exists a set $A \in  \mathcal{R}(f) $ with
\[
A(x) \gg x^{1/3}.
\]
Cilleruelo and Nathanson~\cite{cill-nath06b} recently improved this to 
\[
A(x) \gg x^{\sqrt{2}-1+\varepsilon}
\]
for any $\varepsilon > 0.$

A related problem is the \emph{inverse problem for representation functions}.  Associated to a function $f:\Z \rightarrow \N_0 \cup \{\infty\}$ can be infinitely many sets $A$ of integers such that $A \in \mathcal{R}(f)$.
On the other hand, the semigroup $\N_0$ of nonnegative integers is more rigid than the group \Z\ of integers.  Given $f:\N_0 \rightarrow \N_0$,  there may be a unique set $A \subseteq \N_0$ such that $r_{A,h}(n) = f(n)$ for all $n \in \N_0$.

\bp
Let $f: \N_0 \rightarrow \N_0$ be the representation function of a set of integers.  Determine all sets $A \subseteq \N_0$ such that $r_{A,h}(n) = f(n)$ for all sufficiently large integers $n$.
\ep

The problem was first studied by Nathanson~\cite{nath78}, and subsequently by Lev~\cite{lev04} and others.  There is an excellent survey of additive representation functions by S{\' a}rk{\" o}zy and S{\' o}s~\cite{ss}.

\providecommand{\bysame}{\leavevmode\hbox to3em{\hrulefill}\thinspace}
\providecommand{\MR}{\relax\ifhmode\unskip\space\fi MR }
\providecommand{\MRhref}[2]{%
  \href{http://www.ams.org/mathscinet-getitem?mr=#1}{#2}
}
\providecommand{\href}[2]{#2}


\begin{thebibliography}{1}

\bibitem{cill-nath06b}
J.~Cilleruelo and M.~B.~Nathanson, Dense sets of integers with prescribed representation fucntions, preprint, 2006.


\bibitem{crof67}
H.~T. Croft, \emph{{Research problems, Problem 7, Section 6}}, Mimeographed
  notes, University of Cambridge, 1967.

\bibitem{lev04}
V.~F.~Lev, Reconstructing integer sets from their representation functions,
Electron. J. Combin. \textbf{11} (2004), Research Paper 78, 6 pp. (electronic).



\bibitem{mari69}
J.~Marica, \emph{On a conjecture of {C}onway}, Canad. Math. Bull. \textbf{12}
  (1969), 233--234.



\bibitem{nath78}
M.~B.~Nathanson, Representation functions of sequences in additive number theory, Proc. Amer. Math. Soc. \textbf{72} (1978), 16--20.



\bibitem{nath96}
M.~B.~Nathanson, \emph{Additive Number Theory:  The Classical Bases}, volume 164 of \emph{Graduate Texts in Mathematics}, Springer-Verlag, New York, 1996.

\bibitem{nath00}
M.~B.~Nathanson, \emph{Elementary Methods in Number Theory}, volume 195 of \emph{Graduate Texts in Mathematics}, Springer-Verlag, New York, 1996.


\bibitem{nath03}
M.~B.~Nathanson, Unique representation bases, Acta Arith. \textbf{108} (2003), 1--8.

\bibitem{nath04}
M.~B.~Nathanson, The inverse problem for representation functions of additive bases, in: \emph{Number Theory: New York Seminar 2003}, Springer, New York, 2004, pages 253--262.

\bibitem{nath05}
M.~B.~Nathanson, Every function is the representation function of an additive basis for the integers, Port. Math. (N. S.) \textbf{62} (2005), no. 1, 55--72.

\bibitem{noors06}
M.~B.~Nathanson, K.~O'Bryant, B.~Orosz, I.~Z.~Ruzsa, and M.~Silva, \emph{Binary linear forms}, preprint, 2006.

\bibitem{ss}
A. ~S{\'a}rk{\"o}zy and V.~T.~S{\'o}s, On additive representation functions, in: \emph{The mathematics of Paul Erd\H os, I}, Springer-Verlag, Berlin, 1997, pages 129--150.

\end{thebibliography}
\end{document}